\documentclass[11pt]{article}
\usepackage{amsfonts}
\usepackage{amsmath}

\begin{document}

\begin{center}
\textbf{WINTGEN INEQUALITY FOR STATISTICAL SURFACES}

\bigskip

\textbf{M. Evren AYDIN}$^{1}$\textbf{, Ion MIHAI}$^{2}$

$^{1}$Department of Mathematics, Faculty of Science, Firat University,
Elazig, 23200, Turkey, meaydin@firat.edu.tr

$^{2}$Department of Mathematics, Faculty of Mathematics and Computer
Science, University of Bucharest, Bucharest, 010014, Romania,
imihai@fmi.unibuc.ro\bigskip
\end{center}

\textbf{Abstract. }The Wintgen inequality (1979) is a sharp geometric
inequality for surfaces in the 4-dimensional Euclidean space involving the
Gauss curvature (intrinsic invariant) and the normal curvature and squared
mean curvature (extrinsic invariants), respectively.

In the present paper we obtain a Wintgen inequality for statistical surfaces.
\smallskip



\section{\textbf{Preliminaries}}

For surfaces $M^2$ of the Euclidean space ${\Bbb E}^3$, the Euler inequality
$K\le\Vert H\Vert^2$ is fulfilled, where $K$ is the (intrinsic) Gauss curvature of $M^2$
and $\Vert H\Vert^2$ is the (extrinsic) squared mean curvature of $M^2$.

Furthermore, $K=\Vert H\Vert ^2$
everywhere on $M^2$ if and only if $M^2$ is totally umbilical, or still, by a theorem of Meusnier, if and
only if $M^2$ is (a part of) a plane ${\Bbb E}^2$ or, it is (a part of) a round sphere $S^2$ in ${\Bbb E}^3$.

In 1979, P. Wintgen \cite{W} proved that the Gauss curvature $K$, the squared mean curvature $\left\|H\right\|^2$
and the normal curvature $K^{\bot}$ of any surface $M^{2}$ in ${\Bbb E}^{4}$ always satisfy the inequality
$$K \leq \left\|H\right\|^2-|K^{\bot}|;$$
the equality holds if and only if the ellipse of curvature of $M^{2}$ in ${\Bbb E}^{4}$ is a circle.

The Whitney 2-sphere satisfies the equality case of the Wintgen inequality identically.

A survey containing recent results on surfaces satisfying identically the equality case of Wintgen inequality can be read in 
\cite{Chen-RIGA}. 

Later, the Wintgen inequality was extended by B. Rouxel \cite{R} and by I.V. Guadalupe and L. Rodriguez \cite{GR} independently, for surfaces $M^{2}$ of arbitrary codimension $m$ in real space forms $\widetilde{M}^{2+m}(c)$; namely
$$K \leq \left\|H\right\|^2-|K^{\bot}|+c.$$
The equality case was also investigated.

A corresponding inequality for totally real surfaces in $n$-dimensional complex space forms was obtained in \cite{AM-K}. The equality case was studied and a non-trivial example of a totally real surface satisfying the equality case identically was given (see also \cite{AM-R}).\medskip

In 1999, P.J. De Smet, F. Dillen, L. Verstraelen and L. Vrancken \cite{DDVV} formulated the conjecture on Wintgen inequality for submanifolds of real space forms, which is also known as the {\it DDVV conjecture}. 
   
This conjecture was proven by the authors for submanifolds $M^n$ of arbitrary dimension $n\geq 2$ and codimension $2$
in real space forms $\tilde{M}^{n+2}(c)$ of constant sectional curvature $c$.

Recently, the DDVV conjecture was finally settled for the general case by Z. Lu \cite{L} and independently by J. Ge and Z. Tang \cite{GT}.

One of the present authors obtained generalized Wintgen inequalities for Lagrangian submanifolds in complex space forms \cite{M-N} and Legendrian submanifolds in Sasakian space forms \cite{M-T}, respectively.
 
\section{\textbf{Statistical manifolds and their submanifolds}}

A \textit{statistical manifold} is a Riemannian manifold $\left( \tilde{M}%
^{n+k},\tilde{g}\right) $ of dimension $\left( n+k\right) ,$ endowed with a
pair of torsion-free affine connections $\tilde{\nabla}$ and $\tilde{\nabla}%
^{\ast }$ satisfying%
\begin{equation}
Z\tilde{g}\left( X,Y\right) =\tilde{g}\left( \tilde{\nabla}_{Z}X,Y\right) +%
\tilde{g}\left( X,\tilde{\nabla}_{Z}^{\ast }Y\right) ,  \tag{2.1}
\end{equation}%
for any $X,Y$ and $Z\in \Gamma \left( T\tilde{M}\right) .$ The connections $%
\tilde{\nabla}$ and $\tilde{\nabla}^{\ast }$ are called\textit{\ dual
connections} (see \cite{1}, \cite{NS}), and it is easily shown that $\left( \tilde{\nabla}^{\ast
}\right) ^{\ast }=\tilde{\nabla}.$ The pairing $(\tilde\nabla,\tilde g)$ is said to be
a {\it statistical structure}. If $\left( \tilde{\nabla},\tilde{g}%
\right) $ is a statistical structure on $\tilde{M}^{n+k},$ so is $\left( 
\tilde{\nabla}^{\ast },\tilde{g}\right) $ \cite{1,34}.

On the other hand, any torsion-free affine connection $\tilde{\nabla}$
always has a dual connection given by 
\begin{equation}
\tilde{\nabla}+\tilde{\nabla}^{\ast }=2\tilde{\nabla}^{0},  \tag{2.2}
\end{equation}%
where $\tilde{\nabla}^{0}$ is Levi-Civita connection for $\tilde{M}^{n+k}$.

In affine differential geometry the dual connections are called {\it conjugate connections} (see \cite{Udo}, \cite{Luc}).

Denote by $\tilde{R}$ and $\tilde{R}^{\ast }$ the curvature tensor fields of 
$\tilde{\nabla}$ and $\tilde{\nabla}^{\ast },$ respectively.

A statistical structure $\left( \tilde{\nabla},\tilde{g}\right) $ is said to
be of constant curvature $c\in \mathbb{R}$ if 
\begin{equation}
\tilde{R}\left( X,Y\right) Z=c\left\{ \tilde{g}\left( Y,Z\right) X-\tilde{g}%
\left( X,Z\right) Y\right\} .  \tag{2.3}
\end{equation}%
A statistical structure $\left( \tilde{\nabla},\tilde{g}\right) $ of
constant curvature $0$ is called \textit{a Hessian structure. }

The curvature tensor fields $\tilde{R}$ and $\tilde{R}^{\ast }$ of dual
connections satisfy 
\begin{equation}
\tilde{g}\left( \tilde{R}^{\ast }\left( X,Y\right) Z,W\right) =-\tilde{g}%
\left( Z,\tilde{R}\left( X,Y\right) W\right) .  \tag{2.4}
\end{equation}%
From $\left( 2.4\right) $ it follows immediately that if $\left( \tilde{%
\nabla},\tilde{g}\right) $ is a statistical structure of constant curvature 
$c,$ then $\left( \tilde{\nabla}^{\ast },\tilde{g}\right) $ is also
statistical structure of constant curvature $c$. In particular, if $\left( \tilde{%
\nabla},\tilde{g}\right) $ is Hessian, so is $\left( \tilde{\nabla}^{\ast },%
\tilde{g}\right) $ \cite{13}.

If $\left( \tilde{M}^{n+k},\tilde{g}\right) $ is a statistical manifold and 
$M^n$ a submanifold of dimension $n$ of $\tilde{M}^{n+k},$ then $\left(
M^{n},g\right) $ is also a statistical manifold with the induced connection by $\tilde\nabla$ 
and induced metric $g$. In the case that $\left( \tilde{M}^{n+k},\tilde{g}\right) $ is a
semi-Riemannian manifold, the induced metric $g$ has to be non-degenerate.
For details, see (\cite{33,34}). 

In the geometry of Riemannian submanifolds (see \cite{Chen}), the fundamental equations are
the Gauss and Weingarten formulas and the equations of Gauss, Codazzi and Ricci.

Let denote the set of the sections of the normal bundle to $M^{n}$ by $\Gamma \left( TM^{n\perp }\right).$

In our case, for any $X,Y\in \Gamma \left(TM^{n}\right) ,$ according to  \cite{34}$,$ the corresponding Gauss formulas
are 
$$\tilde{\nabla}_{X}Y =\nabla _{X}Y+h\left( X,Y\right), \eqno (2.5)$$
$$\tilde{\nabla}_{X}^{\ast }Y =\nabla _{X}^{\ast }Y+h^{\ast }\left(
X,Y\right) ,  \eqno (2.6)$$
where $h$, $h^{\ast }:\Gamma (TM^n)\times \Gamma (TM^n)\to \Gamma ({TM^n}^\perp)$ are symmetric and bilinear, called the imbedding
curvature tensor of $M^{n}$ in $\tilde{M}^{n+k}$ for $\tilde{\nabla}$ and
the imbedding curvature tensor of $M^{n}$ in $\tilde{M}^{n+k}$ for $\tilde{%
\nabla}^{\ast },$ respectively.

In \cite{34}, it is also proved that $\left( \nabla ,g\right) $ and $\left(
\nabla ^{\ast },g\right) $ are dual statistical structures on $M^{n}.$

Since $h$ and $h^{\ast }$are bilinear, we have the linear
transformations $A_{\xi }$ and $A_{\xi }^{\ast }$ on $TM^n$ defined by
\begin{equation}
g\left( A_{\xi }X,Y\right) =\tilde{g}\left( h\left( X,Y\right) ,\xi \right), 
\tag{2.7}
\end{equation}%
\begin{equation}
g\left( A_{\xi }^{\ast }X,Y\right) =\tilde{g}\left( h^{\ast }\left(
X,Y\right) ,\xi \right)   \tag{2.8},
\end{equation}%
for any $\xi \in \Gamma \left( TM^{n\perp }\right) $ and $X,Y\in \Gamma
\left( TM^{n}\right) .$ Further, see \cite{34}, the corresponding
Weingarten formulas are
$$\tilde{\nabla}_{X}\xi =-A_{\xi }^{\ast }X+\nabla _{X}^{\perp }\xi , \eqno (2.9)$$
$$\tilde{\nabla}_{X}^{\ast }\xi =-A_{\xi }X+\nabla _{X}^{\ast \perp }\xi,  \eqno (2.10)$$
for any $\xi \in \Gamma \left( TM^{n\perp }\right) $ and $X\in \Gamma \left(
TM^{n}\right) .$ The connections $\nabla _{X}^{\perp }$ and $\nabla
_{X}^{\ast \perp }$ given by $\left( 2.9\right) $ and $\left( 2.10\right) $
are Riemannian dual connections with respect to induced metric on $\Gamma
\left( TM^{n\perp }\right) .$

Let $\left\{ e_{1},...,e_{n}\right\} $ and $\left\{
e_{n+1},...,e_{n+k}\right\} $ be orthonormal tangent and normal frames,
respectively, on $M.$ Then the mean curvature vector fields are defined by%
\begin{equation}
H=\frac{1}{n}\sum_{i=1}^{n}h\left( e_{i},e_{i}\right) =\frac{1}{n}%
\sum_{\alpha =1}^{k}\left( \sum_{i=1}^{n}h_{ii}^{\alpha }\right) e_{n+\alpha
},\text{ }h_{ij}^{\alpha }=\tilde{g}\left( h\left( e_{i},e_{j}\right)
,e_{n+\alpha }\right)  \tag{2.11}
\end{equation}%
and%
\begin{equation}
H^{\ast }=\frac{1}{n}\sum_{i=1}^{n}h^{\ast }\left( e_{i},e_{i}\right) =\frac{%
1}{n}\sum_{\alpha =1}^{k}\left( \sum_{i=1}^{n}h_{ii}^{\ast \alpha }\right)
e_{n+\alpha },\text{ }h_{ij}^{\ast \alpha }=\tilde{g}\left( h^{\ast }\left(
e_{i},e_{j}\right) ,e_{n+\alpha }\right),  \tag{2.2}
\end{equation}%
for $1\leq i,j\leq n$ and $1\le\alpha\le k$ (see also \cite{C-G}).

The corresponding Gauss, Codazzi and Ricci equations are given by the
following result.\medskip

\textbf{Proposition 2.1. }\cite{34}\ \textit{Let }$\tilde{\nabla}$\textit{\
be a dual connection on }$\tilde{M}^{n+k}$\textit{\ and }$\nabla $\textit{\
the induced connection on }$M^{n}.$\textit{\ Let }$\tilde{R}$\textit{\ and }$%
R$\textit{\ be the Riemannian curvature tensors for }$\tilde{\nabla}$
\textit{\ and }$\nabla ,$\textit{\ respectively. Then,}%
$$\tilde{g}\left( \tilde{R}\left( X,Y\right) Z,W\right) =g\left( R\left(
X,Y\right) Z,W\right) +\tilde{g}\left( h\left( X,Z\right) ,h^{\ast }\left(
Y,W\right) \right)-\tilde{g}\left( h^{\ast }\left( X,W\right) ,h\left( Y,Z\right) \right) , 
\eqno (2.13)$$
\begin{eqnarray}
\left( \tilde{R}\left( X,Y\right) Z\right) ^{\perp } &=&\nabla _{X}^{\perp
}h\left( Y,Z\right) -h\left( \nabla _{X}Y,Z\right) -h\left( Y,\nabla
_{X}Z\right)   \notag \\
&&-\left\{ \nabla _{Y}^{\perp }h\left( Y,Z\right) -h\left( \nabla
_{Y}X,Z\right) -h\left( X,\nabla _{Y}Z\right) \right\} ,  \notag
\end{eqnarray}
$$\tilde{g}\left( R^{\perp }\left( X,Y\right) \xi ,\eta \right) =\tilde{g}%
\left( \tilde{R}\left( X,Y\right) \xi ,\eta \right) +g\left( \left[ A_{\xi
}^{\ast },A_{\eta }\right] X,Y\right) ,  \eqno (2.14)$$
\textit{where }$R^{\perp }$\textit{\ is the Riemannian curvature tensor on }$%
TM^{n\perp },$\textit{\ }$\xi ,\eta \in $\textit{\ }$\Gamma \left(
TM^{n\perp }\right) $\textit{\ and }$\left[ A_{\xi }^{\ast },A_{\eta }\right]
=A_{\xi }^{\ast }A_{\eta }-A_{\eta }A_{\xi }^{\ast }.$
\medskip

For the equations of Gauss, Codazzi and Ricci with respect to the dual
connection $\tilde{\nabla}^{\ast }$ on $M^{n}$, we have
\medskip

\textbf{Proposition 2.2. }\textit{Let }$\tilde{\nabla}^{\ast }$\textit{\ be
a dual connection on }$\tilde{M}^{n+k}$\textit{\ and }$\nabla ^{\ast }$%
\textit{\ the induced connection on }$M^{n}.$\textit{\ Let }$\tilde{R}^{\ast
}$\textit{\ and }$R^{\ast }$\textit{\ be the Riemannian curvature tensors
for }$\tilde{\nabla}^{\ast }$\textit{\ and }$\nabla ^{\ast },$\textit{\
respectively. Then,}%
$$\tilde{g}\left( \tilde{R}^{\ast }\left( X,Y\right) Z,W\right)=g\left(
R^{\ast }\left( X,Y\right) Z,W\right) +\tilde{g}\left( h^{\ast }\left(
X,Z\right) ,h\left( Y,W\right) \right)-\tilde{g}\left( h\left( X,W\right) ,h^{\ast }\left( Y,Z\right) \right) , 
\, (2.15)$$
\begin{eqnarray}
\left( \tilde{R}^{\ast }\left( X,Y\right) Z\right) ^{\perp } &=&\nabla
_{X}^{\ast \perp }h^{\ast }\left( Y,Z\right) -h^{\ast }\left( \nabla
_{X}^{\ast }Y,Z\right) -h^{\ast }\left( Y,\nabla _{X}^{\ast }Z\right)  
\notag \\
&&-\left\{ \nabla _{Y}^{\ast \perp }h^{\ast }\left( Y,Z\right) -h^{\ast
}\left( \nabla _{Y}^{\ast }X,Z\right) -h^{\ast }\left( X,\nabla _{Y}^{\ast
}Z\right) \right\}   \notag
\end{eqnarray}%
$$\tilde{g}\left( R^{\ast \perp }\left( X,Y\right) \xi ,\eta \right) =\tilde{g}%
\left( \tilde{R}^{\ast }\left( X,Y\right) \xi ,\eta \right) +g\left( \left[
A_{\xi },A_{\eta }^{\ast }\right] X,Y\right) ,  \eqno (2.16)$$
\textit{where }$R^{\ast \perp }$\textit{\ is the Riemannian curvature tensor for 
}$\nabla ^{\perp \ast }$\textit{\ on }$TM^{n\perp },$\textit{\ }$\xi ,\eta
\in $\textit{\ }$\Gamma \left( TM^{n\perp }\right) $\textit{\ and }$\left[
A_{\xi },A_{\eta }^{\ast }\right] =A_{\xi }A_{\eta }^{\ast }-A_{\eta }^{\ast
}A_{\xi }.$\medskip

Geometric inequalities for statistical submanifolds in statistical manifolds with constant curvature were obtained
in \cite{2}.

\section{\textbf{Sectional curvature for statistical manifolds}}

Let $\left( M^{n},g\right) $ be a statistical manifold of dimension $n$
endowed with dual connections $\tilde{\nabla}$ and $\tilde{\nabla}^{\ast }.$
Unfortunately, the $\left( 0,4\right)$-tensor field $g\left( R\left(
X,Y\right) Z,W\right) $ is not skew-symmetric relative to $Z$ and $W.$ Then
we cannot define a sectional curvature on $M^{n}$ by the standard definition.

We shall define a skew-symmetric $\left( 0,4\right)$-tensor field on $M^{n}$ by
\begin{equation*}
T\left( X,Y,Z,W\right) =\frac{1}{2}\left[ g\left( R\left( X,Y\right)Z,W\right)+g\left( R^{\ast }\left( X,Y\right) Z,W\right) \right], 
\end{equation*}%
for all $X,Y,Z,W\in \Gamma \left( TM^{n}\right) .$

Then we are able to define a sectional curvature on $M^{n}$ by the formula%
\begin{equation*}
K\left( X\wedge Y\right) =\frac{T\left( X,Y,X,Y\right) }{g\left( X,X\right)
g\left( Y,Y\right) -g^{2}\left( X,Y\right) },
\end{equation*}%
for any linearly independent tangent vectors $X,Y$ at $p\in M^{n}.$

We want to point-out that this definition has the opposite sign that the
sectional curvature defined by B. Opozda \cite{Opozda8}. Another sectional curvature
was considered in \cite{Opozda} (see also \cite{S}).

In particular, for a statistical surface $M^{2}$, we can define a Gauss
curvature by
\begin{equation*}
G=K\left( e_{1}\wedge e_{2}\right),
\end{equation*}%
for any orthonormal frame $\left\{ e_{1},e_{2}\right\} $ on $M^{2}.$

Analogously, we shall consider a normal curvature of a statistical surface $%
M^{2}$ in an orientable $4$-dimensional statistical manifold $\tilde M^{4}.$ Let $%
\left\{ e_{1},e_{2},e_{3},e_{4}\right\} $ be a positive oriented orthonormal
frame on $\tilde{M}^{4}$ such that $e_{1},e_{2}$ are tangent to $M^{2}.$
Let 
\begin{equation*}
G^{\perp }=\frac{1}{2}\left[ g\left( R^{\perp }\left( e_{1},e_{2}\right)
e_{3},e_{4}\right)+g\left( R^{\ast \perp }\left( e_{1},e_{2}\right)
e_{3},e_{4}\right) \right]
\end{equation*}
be a normal curvature of $M^{2}.$

Remark that $\left\vert G^{\perp }\right\vert $ does not depend on the
orientation of the statistical manifold. Then $\left\vert G^{\perp
}\right\vert $ can be defined for any surface $M^{2}$ of any $4$-dimensional
statistical manifold.\bigskip

We state a version of Euler inequality for surfaces in 3-dimensional statistical manifolds of constant curvature.\medskip

{\bf Proposition 3.1.} {\it Let $M^2$ be surface in a $3$-dimensional statistical manifold of constant curvature $c$. Then its Gauss curvature satisfies}:
$$G\le 2||H||\cdot ||H^*||-c.$$

{\it Proof.} Let $p\in M^2$ and $e_3$ be a unit normal vector to $M^2$ at $p$. We can choose an orthonormal basis $\{e_1,e_2\}$ of $T_pM^2$ such that $h^0(e_1,e_2)=0$, where $h^0$ is the second fundamental form of $M^2$ (with respect to the Levi-Civita connection). Then $h_{12}^{3}+h_{12}^{*3}=0$.

The Gauss equations for $\nabla$ and $\nabla^*$ imply 
$$G=-c-\frac 12(h_{11}^3h_{22}^{*3}+h_{11}^{*3}h_{22}^3)+h_{12}^3h_{12}^{*3}.$$

Applying the Cauchy-Schwarz inequality, it follows that
$$G\le -c+\frac 12\sqrt{(h_{11}^3+h_{22}^3)^2(h_{11}^{*3}+h_{22}^{*3})^2}-(h_{12}^3)^2.$$

But in our case $4||H||^2=(h_{11}^3+h_{22}^3)^2$ and $4||H^*||^2=(h_{11}^{*3}+h_{22}^{*3})^2$. Therefore
$$G\le -c+2||H||\cdot ||H^*||.$$

\textbf{Example 1.} (A trivial example) Recall Lemma 5.3 of Furuhata \cite{13}.

Let $\left( \mathbb{H},\tilde{\nabla},\tilde{g}\right) $ be a Hessian
manifold of constant Hessian curvature $\tilde{c}\neq 0,$ $\left( M,\nabla ,g\right) $ a 
trivial Hessian manifold and $f:M\longrightarrow \mathbb{H}$ a
statistical immmersion of codimension one. Then one has:
\begin{equation*}
A^{\ast }=0,\text{ \ \ }h^{\ast }=0,\text{ \ \ }\left\Vert H^{\ast
}\right\Vert =0.
\end{equation*}

Thus, if dim $M=2$, the immersion $f$ of codimension one satisfies the equality case of the
statistical version of Euler inequality given by Proposition 3.1.\medskip 

\textbf{Example 2.} Let $\left( \mathbb{H}^{3},\tilde{g}\right) $ be the
upper half space of constant sectional curvature $-1$, i.e.,
$${\Bbb H}^3 =\left\{ y=\left( y^{1},y^{2},y^{3}\right) \in \mathbb{R}^{3}:y^{3}>0\right\}, \quad
\tilde{g}=\left( y^{3}\right) ^{-2}\sum_{k=1}^3 dy^{k}dy^{k}.$$

An affine connection $\tilde{\nabla}$ on $\mathbb{H}$ is given by%
$$\tilde{\nabla}_{\frac{\partial }{\partial y^{3}}}\frac{\partial }{\partial
y^{3}}=\left( y^{3}\right) ^{-1}\frac{\partial }{\partial y^{3}}, \quad
\tilde{\nabla}_{\frac{\partial }{\partial y^{i}}}\frac{\partial }{\partial
y^{j}}=2\delta _{ij}\left( y^{3}\right) ^{-1}\frac{\partial }{\partial y^{3}}, \quad
\tilde{\nabla}_{\frac{\partial }{\partial y^{i}}}\frac{\partial }{\partial
y^{3}}=\tilde{\nabla}_{\frac{\partial }{\partial y^{3}}}\frac{\partial }{%
\partial y^{j}}=0,$$
where $i,j=1,2.$ The curvature tensor field $\tilde{R}$ of $\tilde{\nabla}$
is identically zero, i.e., $c=0$. Thus $\left( \mathbb{H}^{3},\tilde{\nabla},\tilde{g}\right) $ 
is a Hessian manifold of constant Hessian curvature 4.

Now let consider a horosphere $M^{2}$ in $\mathbb{H}^{3}$ having null
Gauss curvature, i.e., $G\equiv 0.$ (For details, see \cite{Lopez}).
If $f:M^{2}\longrightarrow \mathbb{H}^{3}$ is a statistical immersion of
codimension one, then, by using Lemma 4.1 of \cite{Min}, we deduce $A^{\ast }=0$,
and then $H^{\ast }=0.$ This implies that the horosphere $M^{2}$ satisfies the
equality case of the statistical version of Euler inequality given by
Proposition 3.1.

\section{Wintgen inequality for statistical surfaces in a 4-dimensional
statistical manifold of constant curvature}

Let $\left( \tilde{M}^{4},c\right) $ be a statistical manifold of constant
curvature $c$ and $M^{2}$ a statistical surface in $\left( \tilde{M}%
^{4},c\right) .$

We shall prove a Wintgen inequality for the surfaces $M^{2}$ in $\left( 
\tilde{M}^{4},c\right) .$ The Gauss curvature $G$ of $M^{2}$ is given by%
\begin{equation*}
G=\frac{1}{2}\left[ g\left( R\left( e_{1},e_{2}\right) e_{1},e_{2}\right)
+g\left( R^{\ast }\left( e_{1},e_{2}\right) e_{1},e_{2}\right) \right] .
\end{equation*}%
By the Gauss equation we have 
\begin{equation*}
g\left( R\left( e_{1},e_{2}\right) e_{1},e_{2}\right) =g\left( \tilde{R}%
\left( e_{1},e_{2}\right) e_{1},e_{2}\right) -g\left( h\left(
e_{1},e_{1}\right) ,h^{\ast }\left( e_{2},e_{2}\right) \right) +g\left(
h^{\ast }\left( e_{1},e_{2}\right) ,h\left( e_{1},e_{2}\right) \right) ,
\end{equation*}%
or equivalently,%
\begin{equation*}
g\left( R\left( e_{1},e_{2}\right) e_{1},e_{2}\right)
=-c-h_{11}^{3}h_{22}^{\ast 3}-h_{11}^{4}h_{22}^{\ast 4}+h_{12}^{\ast
3}h_{12}^{3}+h_{12}^{\ast 4}h_{12}^{4}.
\end{equation*}
Analogously
\begin{equation*}
g\left( R^{\ast }\left( e_{1},e_{2}\right) e_{1},e_{2}\right)
=-c-h_{11}^{\ast 3}h_{22}^{3}-h_{11}^{\ast
4}h_{22}^{4}+h_{12}^{3}h_{12}^{\ast 3}+h_{12}^{4}h_{12}^{\ast 4}.
\end{equation*}%
It follows that%
\begin{equation*}
G=-c-\frac{1}{2}\left[ h_{11}^{3}h_{22}^{\ast 3}+h_{11}^{\ast
3}h_{22}^{3}+h_{11}^{4}h_{22}^{\ast 4}+h_{11}^{\ast 4}h_{22}^{4}\right]
+h_{12}^{3}h_{12}^{\ast 3}+h_{12}^{4}h_{12}^{\ast 4}.
\end{equation*}%
The normal curvature $G^{\perp }$ of $M^{2}$ is given by%
\begin{equation*}
2\left\vert G^{\perp }\right\vert =\frac{1}{2}\left\vert g\left( R^{\perp
}\left( e_{1},e_{2}\right) e_{3},e_{4}\right) +g\left( R^{\ast \perp }\left(
e_{1},e_{2}\right) e_{3},e_{4}\right) \right\vert .
\end{equation*}%
By the Ricci equations with respect to $\tilde{\nabla}$ and $\tilde{\nabla},$
respectively, we get%
$$2\left\vert G^{\perp }\right\vert =|g([A_{e_3}^*,A_{e_4}]e_1,e_2)+g([A_{e_3},A_{e_4}^*]e_1,e_2)|$$
$$=\left\vert h_{12}^{\ast 3}\left(h_{11}^{4}-h_{22}^{4}\right) -h_{12}^{\ast 4}\left(h_{11}^{3}-h_{22}^{3}\right)+h_{12}^{3}
\left( h_{11}^{\ast 4}-h_{22}^{\ast 4}\right)-h_{12}^{4}\left( h_{11}^{\ast 3}-h_{22}^{\ast 3}\right)\right\vert .$$

In order to estimate $\left\vert G^{\perp }\right\vert ,$ we shall use the
inequalities
\begin{equation*}
\pm 4ab\leq a^{2}+4b^{2},\text{ }a,b\in \mathbb{R}.
\end{equation*}%
Then we have
$$2\left\vert G^{\perp }\right\vert \leq\frac 14[(h_{11}^3-h_{22}^3)^2+(h_{11}^4-h_{22}^4)^2+(h_{11}^{*3}-h_{22}^{*3})^2
+(h_{11}^{*4}-h_{22}^{*4})^2]+(h_{12}^3)^2+(h_{12}^4)^2+(h_{12}^{*3})^2+(h_{12}^{*4})^2$$
$$=\frac{1}{4}\left[ \left\Vert h_{11}-h_{22}\right\Vert ^{2}+\left\Vert h_{11}^{\ast }
-h_{22}^{\ast }\right\Vert ^{2}\right] +\left\Vert h_{12}\right\Vert ^{2}+\left\Vert
h_{12}^{\ast }\right\Vert ^{2},$$
which yields that
\begin{eqnarray*}
2\left\vert G^{\perp }\right\vert &\leq &\frac{1}{4}\left[ \left\Vert
h_{11}+h_{22}\right\Vert ^{2}+\left\Vert h_{11}^{\ast }+h_{22}^{\ast
}\right\Vert ^{2}\right] -g\left( h_{11},h_{22}\right) -g\left( h_{11}^{\ast
},h_{22}^{\ast }\right) +\left\Vert h_{12}\right\Vert ^{2}+\left\Vert
h_{12}^{\ast }\right\Vert ^{2} \\
&=&\left\Vert H\right\Vert ^{2}+\left\Vert H^{\ast }\right\Vert ^{2}-\left(
h_{11}^{3}+h_{11}^{\ast 3}\right) \left( h_{22}^{3}+h_{22}^{\ast 3}\right)
-\left( h_{11}^{4}+h_{11}^{\ast 4}\right) \left( h_{22}^{4}+h_{22}^{\ast
4}\right)  \\
&&+h_{11}^{3}h_{22}^{\ast 3}+h_{11}^{\ast
3}h_{22}^{3}+h_{11}^{4}h_{22}^{\ast 4}+h_{11}^{\ast 4}h_{22}^{4}+\left( h_{12}^{3}\right) ^{2}+\left( h_{12}^{4}\right) ^{2}+
\left(h_{12}^{\ast 3}\right) ^{2}+\left( h_{12}^{\ast 4}\right) ^{2}.
\end{eqnarray*}

It is known that $2h^{0}=h+h^{\ast },$ where $h^{0}$ denotes the second
fundamental form of $M^{2}$ with respect to the Levi-Civita connection $%
\tilde{\nabla}^{0}$ on $\left( \tilde{M}^{4},c\right) .$

Then we can write%
\begin{eqnarray*}
2\left\vert G^{\perp }\right\vert &\leq &\left\Vert H\right\Vert
^{2}+\left\Vert H^{\ast }\right\Vert ^{2}-4\left(
h_{11}^{03}h_{22}^{03}+h_{11}^{04}h_{22}^{04}\right) -2G-2c \\
&&+2h_{12}^{3}h_{12}^{\ast 3}+2h_{12}^{4}h_{12}^{\ast 4}+\left(
h_{12}^{3}\right) ^{2}+\left( h_{12}^{\ast 3}\right) ^{2}+\left(
h_{12}^{4}\right) ^{2}+\left( h_{12}^{\ast 4}\right) ^{2}.
\end{eqnarray*}%
Recall the Gauss equation for the Levi-Civita connection 
\begin{equation*}
\tilde{K}^{0}\left( e_{1}\wedge e_{2}\right)
=G^{0}-h_{11}^{03}h_{22}^{03}-h_{11}^{04}h_{22}^{04}+\left(
h_{12}^{03}\right) ^{2}+\left( h_{12}^{04}\right) ^{2},
\end{equation*}%
where $\tilde{K}^{0}\left( e_{1}\wedge e_{2}\right) $ is the sectional
curvature of $M^{2}$ in $\left( \tilde{M}^{4},\tilde{\nabla}^{0}\right) $
and $G^{0}$ its Gaussian curvature with respect to the Levi-Civita
connection.

Consequently we have%
\begin{eqnarray*}
2\left\vert G^{\perp }\right\vert &\leq &\left\Vert H\right\Vert
^{2}+\left\Vert H^{\ast }\right\Vert ^{2}-4G^{0}+4\tilde{K}^{0}\left(
e_{1}\wedge e_{2}\right) -2G-2c- \\
&&-\left( h_{12}^{3}+h_{12}^{\ast 3}\right) ^{2}-\left(
h_{12}^{4}+h_{12}^{\ast 4}\right) ^{2}+2h_{12}^{3}h_{12}^{\ast
3}+2h_{12}^{4}h_{12}^{\ast 4} \\
&&+\left( h_{12}^{3}\right) ^{2}+\left( h_{12}^{\ast 3}\right) ^{2}+\left(
h_{12}^{4}\right) ^{2}+\left( h_{12}^{\ast 4}\right) ^{2} \\
&=&\left\Vert H\right\Vert ^{2}+\left\Vert H^{\ast }\right\Vert ^{2}-4G^{0}+4%
\tilde{K}^{0}\left( e_{1}\wedge e_{2}\right) -2G-2c.
\end{eqnarray*}%
Summing up, we state the following.\medskip

\textbf{Theorem 4.1.} \textit{Let }$M^{2}$ \textit{be a statistical surface in a} 
$4$\textit{-dimensional statistical manifold} $\left( \tilde{M}^{4},c\right)$ 
\textit{of constant curvature }$c.$\textit{\ Then} 
\begin{equation*}
G+\left\vert G^{\perp }\right\vert +2G^{0}\leq \frac{1}{2}\left( \left\Vert
H\right\Vert ^{2}+\left\Vert H^{\ast }\right\Vert ^{2}\right) -c+2\tilde{K}%
^{0}\left( e_{1}\wedge e_{2}\right) .
\end{equation*}

In particular, for $c=0$ we derive the following.\medskip

\textbf{Corollary 4.2. }\textit{Let }$M^{2}$ \textit{be a statistical surface of
a Hessian }$4$\textit{-dimensional statistical manifold }$\tilde{M}^{4}$ {\it 
of Hessian curvature $0$. Then}:
\begin{equation*}
G+\left\vert G^{\perp }\right\vert +2G^{0}\leq \frac{1}{2}\left( \left\Vert
H\right\Vert ^{2}+\left\Vert H^{\ast }\right\Vert ^{2}\right) .
\end{equation*}

{\bf Acknowledgements.} This paper was written during the visit of the second author to Firat University,
Turkey, in April-May 2015. The second author was supported by the Scientific
and Technical Research Council of Turkey (TUBITAK) for Advanced Fellowships ¨
Programme.

\end{document}